\documentclass[reqno,10pt,centertags,draft]{amsart}
\usepackage{amsmath,amsthm,amscd,amssymb,latexsym,upref}

\newcommand{\bbN}{{\mathbb{N}}}
\newcommand{\bbR}{{\mathbb{R}}}

\newcommand{\cA}{\mathcal {A}}
\newcommand{\cB}{{\mathcal B}}
\newcommand{\cC}{{\mathcal C}}
\newcommand{\cD}{{\mathcal D}}

\newcommand{\cH}{{\mathcal H}}

\newcommand{\cS}{{\mathcal S}}

\newcommand{\lb}{\label}
\newcommand{\f}{\frac}

\newcommand{\ol}{\overline}

\newcommand{\wti}{\widetilde}

\newcommand{\loc}{\text{\rm{loc}}}

\newcommand{\ran}{\text{\rm{ran}}}

\newcommand{\dom}{\text{\rm{dom}}}
\newcommand{\ess}{\text{\rm{ess}}}
\newcommand{\Ell}{\text{\rm{Ell}}}
\newcommand{\Ext}{\text{\rm{Ext}}}
\newcommand{\ac}{\text{\rm{ac}}}

\newcommand{\bi}{\bibitem}
\newcommand{\hatt}{\widehat}
\newcommand{\beq}{\begin{equation}}
\newcommand{\eeq}{\end{equation}}
\newcommand{\ba}{\begin{align}}
\newcommand{\ea}{\end{align}}

\newcommand{\ord}{\text{\rm{ord}}}

\renewcommand{\ge}{\geqslant}
\renewcommand{\le}{\leqslant}

\newcommand{\Om}{\Omega}

\newcommand{\N}{\mathbb{N}}

\allowdisplaybreaks
\numberwithin{equation}{section}

\newcommand{\gH}{\mathfrak{H}}

\allowdisplaybreaks
\numberwithin{equation}{section}

\newtheorem{theorem}{Theorem}[section]

\newtheorem{proposition}[theorem]{Proposition}
\newtheorem{corollary}[theorem]{Corollary}
\newtheorem{hypothesis}[theorem]{Hypothesis}
\theoremstyle{definition}
\newtheorem{definition}[theorem]{Definition}
\newtheorem{remark}[theorem]{Remark}

\begin{document}

\title[Spectral Theory of Elliptic Operators in Exterior Domains] {Spectral Theory 
of Elliptic Operators in Exterior Domains}
\author[F.\ Gesztesy and M.\ M.\ Malamud]{Fritz Gesztesy and Mark M.\ Malamud}
\address{Department of Mathematics,
University of Missouri, Columbia, MO 65211, USA}
\email{fritz@math.missouri.edu}
\urladdr{http://www.math.missouri.edu/personnel/faculty/gesztesyf.html}
\address{Mathematics, Institute of Applied Mathematics and Mechanics,
R, Luxemburg str. 74, Donetsk 83114, Ukraine}
\email{mmm@telenet.dn.ua}

\begin{abstract}
We consider various closed (and self-adjoint) extensions of elliptic differential expressions of the type $\cA=\sum_{0\le |\alpha|,|\beta|\le m}(-1)^\alpha
D^\alpha a_{\alpha, \beta}(x)D^\beta$,
$a_{\alpha, \beta}(\cdot)\in C^{\infty}({\overline\Omega})$, on smooth (bounded or unbounded) domains $\Om$ in $\bbR^n$ with compact boundary $\partial\Omega$.
Using the concept of boundary triples and operator-valued Weyl--Titchmarsh functions, we prove various trace ideal properties of powers of resolvent differences of these closed
realizations of $\cA$ and derive estimates on eigenvalues of certain self-adjoint realizations in spectral gaps of the Dirichlet realization.

Our results extend classical theorems due to Vi{\u s}ik, Povzner, Birman, and Grubb.
\end{abstract}

\maketitle

\section{Introduction}\label{s1}

Let $\Omega$ be an open domain in $\bbR^n$ (bounded or unbounded) with compact boundary $\partial\Omega$. Throughout we assume that $\partial\Omega$ is an  
$(n-1)$-dimensional (not necessarily connected) $C^{\infty}$-manifold. Let $\cA$
be the differential expression
\begin{equation}\label{2.10}
\cA=\sum_{0\le |\alpha|,|\beta|\le m}(-1)^\alpha D^\alpha
a_{\alpha, \beta}(x)D^\beta, \quad a_{\alpha, \beta}(\cdot)\in
C^{\infty}({\overline\Omega}),
\end{equation}
 $\ord(\cA) = 2m$,  which is  elliptic in ${\overline\Omega}$.
Moreover, we assume that  $\cA$ is properly elliptic in $\ol \Om$
(which is automatically satisfied if either $n>2$ or the symbol of
$\cA$ is real, cf.\ \cite{LM72}). In addition to \eqref{2.10} we consider its formal adjoint
$\cA^{\top} = \sum_{0\le |\alpha|,|\beta|\le m}(-1)^\alpha D^\beta
{\overline{a_{\beta,\alpha}(x)}} D^\beta$, which is also properly
elliptic in ${\overline\Omega}$ (cf.\ \cite{LM72}).

Denote by $A=A_{\min}$ $(A^{\top} = A^{\top}_{\min})$ the minimal operator
associated in $L^2(\Omega)$ with the differential expression
$\cA$ (resp., $\cA^{\top}$), that is, the closure of $\cA$ defined on
$C_0^{\infty}(\Omega)$.  The maximal operators $A_{\max}$ and
$A^{\top}_{\max}$ are then defined by $A_{\max}=(A^{\top}_{\min})^* = (A^{\top})^*$ and
$A^{\top}_{\max}=(A_{\min})^*=A^*$, respectively.
We emphasize that $H^{2m}(\Omega)\subset \dom(A_{\max})\subset
H^{2m}_{\loc}(\Omega)$, while $\dom(A_{\max})\not =H^{2m}(\Omega)$.

After the pioneering work by Vishik \cite{Vi63}, nonlocal boundary
value problems of the form $A_{\max}u = f$,
$(\partial u/\partial n - Ku)\upharpoonright\partial\Omega =0$ for elliptic operators \eqref{2.10}
(with $m=1$) in bounded domains were considered by numerous
authors (see, e.g., \cite{Be65, Gr68} and the references therein). Vishik
was the first to consider these  problems
in the framework of extension theory of dual pairs of operators.
Starting with a formula for the domain $\dom(A_{\max})$ of
$A_{\max}$, he applied it to an appropriate regularization of the
classical Green's  formula, using the Calderon operator. The
latter allowed him  to extend the Green's  formula  from
$H^{2m}(\Omega)$ to $\dom(A_{\max})$. The next fundamental
contribution to the subject was made by Grubb \cite{Gr68}.
Using the theory of Lions and Magenes \cite{LM72}, Grubb
substantially extended and completed the results of \cite{Vi63}. In particular, Grubb
obtained the (regularized) Green's  formula which (in the special case $m=1$)
reads as follows:
$$
(A_{\max}u, v)_{L^2(\Om)} - (u, A^{\top}_{\max}v)_{L^2(\Om)}
 = \big({\widetilde\Gamma}_{\Om,1} u,{\widetilde\Gamma}^{\top}_{\Om,0}
 v\big)_{1/2,-1/2}
 - \big({\widetilde\Gamma}_{\Om,0} u,{\widetilde\Gamma}^{\top} _{\Om,1}
 v\big)_{-1/2,1/2}.
$$
Here $(\cdot, \cdot)_{s,-s}$ denotes the duality pairing between
$H^{s}(\partial \Omega)$ and $H^{-s}(\partial \Omega)$, $u\in
\dom(A_{\max})$, $v\in \dom(A^{\top}_{\max})$, and $\wti \Gamma_{\Om,0}$,
$\wti \Gamma_{\Om,1}$, ${\Gamma}^{\top}_{\Om,0}$,
and ${\widetilde\Gamma}^{\top} _{\Om,1}$ are regularized trace
operators, having the properties
$$
{\widetilde\Gamma}_{\Om,1} \colon \cD(A_{\max}) \to
H^{1/2}(\partial\Omega), \quad
{\widetilde\Gamma^{\top}}_{\Om,0} \colon \cD(A_{\max})
\to H^{-1/2}(\partial\Omega), \quad  \ran \big(
\big({\widetilde\Gamma}_{\Om,1},{\widetilde\Gamma^{\top}}_{\Om,0}\big)\big) =
H^{1/2}(\partial \Omega)\times H^{-1/2}(\partial \Omega).
$$
Later, we will use a somewhat different approach (cf.\ Proposition \ref{prop3.2}).

On the other hand, during the past three decades a new approach to
the extension theory, based on the concept of a boundary triple and
the corresponding operator-valued Weyl--Titchmarsh function, was
developed in \cite{DM91} (cf.\ the references therein for the
symmetric case) and in \cite{MM02} (in the case of dual pairs). In
this paper we apply some results and technique from \cite{DM91}
and \cite{MM02} to elliptic operators on unbounded domains. The
most important ingredients from the elliptic theory we need are
the regularized Green's  formula and {\it a priori}
coercivity-type estimates for the elliptic realizations ${\hatt
A}_B$ of $\cA$ (see \eqref{2.54} below). To obtain the latter on
unbounded domains one needs additional restrictions on the
coefficients of $\cA$, since an elliptic realization is not
necessarily coercive. Here we restrict ourselves to the case of
bounded coefficients $a_{\alpha, \beta}(\cdot)$. Using the
formalism of boundary triples and the corresponding
operator-valued Weyl--Titchmarsh functions in \cite {DM91, MM02},
we investigate the resolvent difference of two realizations and
complement the results of Povzner \cite{Po67}, Birman \cite{Bi62},
and Grubb \cite{Gr83} in this direction. In addition, assuming
$A_{\min}>0$, we compute the number of negative eigenvalues of a
realization $A_K $ and the number of eigenvalues of $A_K $ within
spectral gaps of the Dirichlet realization ${\hatt A}_{\gamma_D}$, 
where $\gamma_D =\{\gamma_j\}_0^{m-1}$.

{\bf Notations.}
$\gH$ and $\cH$ represent complex, separable Hilbert spaces;
$\cB(\cH)$, $\cB_\infty(\cH)$, and $\cC(\cH)$ denote the sets of bounded, compact, and closed linear operators in $\cH$; $\dom(\cdot)$, $\ran(\cdot)$, and $\ker(\cdot)$ denote the domain, range, and kernel of a linear operator, $\rho(\cdot)$ and $\sigma(\cdot)$ stand for the resolvent
set and spectrum of a linear operator.
As usual, $C^{\infty}(\Omega)$ denotes the set of infinity
differentiable functions in the domain $\Omega$,
$C^{\infty}_0(\Omega)$ the subset of $C^{\infty}(\Omega)$-functions of compact support in $\Omega$; $C_b({\Omega})=C({\Omega})\cap L^{\infty}({\Omega})$,
$C_{u}(\Omega)$ the set of uniformly continuous functions  in $\Omega$,
$C_{ub}(\Omega) = C_{u}(\Omega) \cap C_{b}(\Omega)$, and $H^s(\Omega)$ the usual Sobolev spaces.

\section{Dual pairs, boundary triples, and operator-valued Weyl--Titchmarsh functions}
\lb{s2}

\begin{definition} \lb{d3.1}
Let  $A$ and $A^{\top}$ be densely defined (not necessarily
closed) linear operators in  $\gH$. Then $A$ and $A^{\top}$ form a
{\it dual pair $($DP\,$)$} $\{A,A^{\top}\}$ in $\gH$ if $(Af,g) =
(f, A^{\top}g)$ for all $f\in \dom (A)$, $g\in \dom (A^{\top})$.
An operator $\widetilde A$ is called a {\it proper extension of
the DP $\{A, A^{\top}\}$}, and we write $\widetilde A\in Ext\{A, A^{\top}\}$, if
$A\subsetneqq \widetilde A\subsetneqq (A^{\top})^*$.
\end{definition}

\begin{definition} (cf.\ \cite{LS83}, \cite{MM02}) \lb{d3.4}
$(i)$ Let $\gH$, $\cH_0$, and $\cH_1$ be complex, separable Hilbert spaces and
$$
\Gamma^{\top}=\begin{pmatrix}\Gamma^{\top}_0\cr\Gamma^{\top}_1\endpmatrix \colon
\dom ((A^{\top})^*)\to \cH_0\oplus
\cH_1 \, \text{ and } \, \Gamma = \pmatrix\Gamma_0\cr\Gamma_1\end{pmatrix}
 \colon \dom (A^*)\to \cH_1\oplus \cH_0
$$
be linear mappings. Then $\Pi= \{\cH_0\oplus \cH_1, \Gamma^{\top}, \Gamma \}$  is
called a {\it boundary triple} for the dual pair  $\{A,A^{\top}\}$ if $\Gamma^{\top}$
and $\Gamma$ are surjective and the Green's  identity holds,
$$
\big((A^{\top})^*f,g\big)_{\gH}-(f,A^*g)_{\gH} = \big(\Gamma^{\top}_1 f, \Gamma_0 g\big)_{\cH_1}
- \big(\Gamma^{\top}_0 f,\Gamma_1 g\big)_{\cH_0}, \quad   f\in \dom\big((A^{\top})^*\big),\,
g\in \dom(A^*).
$$
We set $A_0 = (A^{\top})^*\upharpoonright \ker (\Gamma^{\top}_0)$ and
$A^{\top}_0 = A^*\upharpoonright \ker (\Gamma_0)$. \\
$(ii)$ The operator-valued function $M_{\Pi} (z )$ defined by
$$
\Gamma^{\top}_1 f_{z }=M_{\Pi}(z )\Gamma^{\top}_0 f_{z }, \quad
f_{z } \in \ker \big((A^{\top})^*-z\big), \;\;  z \in \rho (A_0),
$$
is called the {\it Weyl--Titchmarsh  function}
corresponding to the boundary triple $\Pi$.
\end{definition}

Due to Green's  identity, $(A_{\min}u,v)_{L^2(\Om)}= (u, A^{\top}_{\min}v)_{L^2(\Om)}$,
$u, v\in C_0^{\infty}(\Omega)$, the operators $A$ and $A^{\top}$ form a
dual pair $\{A,A^{\top}\}$ of elliptic operators in $L^2(\Om)$.
Any proper extension $\widetilde A\in \Ext\{A,A^{\top}\}$ of $\{A,A^{\top}\}$ is
called a realization of $\cA$. Clearly, any realization  $\widetilde A$ of $\cA$ is
closable. We equip $\dom(A_{\max})$ and $\dom(A^{\top}_{\max})$ with the corresponding graph norms. It is known (cf.\ \cite{Be68, LM72}) that if a  domain
$\Omega$ is bounded, then $ \dom(A_{\min}) = \dom(A^{\top}_{\min}) =
H^{2m}_0(\Omega)$,  where the norms in $\dom(A_{\min})$ and
$H^{2m}_0(\Omega)$ are equivalent. Denote by $\gamma_j$ the mappings
$\gamma_j\colon  C^{\infty}({\overline{\Omega}})\to
C^{\infty}(\partial\Omega)$, $\gamma_j u = \gamma_0
({\partial^j u}/{\partial n^j})= {\partial^j u}/{\partial n^j}
\upharpoonright {\partial\Omega}$, $1\le j \le m-1$, $\gamma_0 u=u \upharpoonright {\partial\Omega}$,
where $n$ stands for the interior normal to $\partial\Omega$.
Next we introduce the boundary operators $B_j$ as
\begin{equation}\label{2.17}
B_j u=\sum_{0 \le |\beta|\le m_j}b_{j\beta}\gamma_0(D^\beta u), \quad
b_{j\beta}(\cdot)\in C^{\infty}(\partial\Omega), \quad \ord (B_j) = m_j\le 2m-1.
\end{equation}
Here $B_j: C^{\infty}({\overline{\Omega}})\to C^{\infty}(\partial\Omega)$ will eventually be extended to appropriate Sobolev spaces $H^{s}(\Omega)$ and in some cases to
$\cD(A_{\max})$. $B_j$ in \eqref{2.17} can also be rewritten as $B_ju = b_j\gamma_{m_j}u\  +\  \sum_{0\le k\le m_j-1}
T_{j,k}\gamma_k u$, where $b_j(\cdot)\in C^{\infty}(\partial\Omega)$ and $T_{j,k}$ are
tangential differential operators in $\partial\Omega$ of orders
$\ord (T_{j,k}) \le m_j-k$ with $C^{\infty}(\partial\Omega)$-coefficients.

With any elliptic operator $A$ \eqref{2.10} and a system $B =\{B_j\}_{j=1}^{m-1}$
we associate the operator ${\hatt A}_B$   defined by
\begin{equation}\label{2.54}
{\hatt A}_B= A_{\max}\upharpoonright \dom({\hatt A}_B), \quad
\dom({\hatt A}_B)=H_B^{2m}(\Omega)= \{u\in H^{2m}(\Omega) \,|\, Bu =0\}.
\end{equation}
Our considerations are based on \cite[Thm.\ 2.2.1]{LM72}.
According to this result, for any elliptic differential expression
$\cA$ in \eqref{2.10} and any normal system $\{B_j\}^{m-1}_{j=0}$ on
$\partial\Omega$ given by \eqref{2.17}, there exists a  system of
boundary operators $\{C_j\}^{m-1}_{j=0}$,  $\ord(C_j) = \mu_j\le 2m-1$, such that the system
$\{B_0,\ldots,B_{m-1},C_0,\ldots,C_{m-1}\}$ is a Dirichlet system
of order $2m$  and  another Dirichlet system of boundary operators $\{B^{\top}_j\}^{m-1}_{j=0}\cup
\{C^{\top}_j\}^{m-1}_{j=0}$, such that the following Green's formula hold
\begin{equation}\label{2.19}
\big(\cA u,  v\big)_{L^2(\Omega)} - \big(u, \cA
^{\top}v\big)_{L^2(\Omega)} = \sum_{0\le l\le m-1} \big[\big(C_j u,
B^{\top}_j v\big)_{L^2(\partial\Omega)}  - \big( B_j u, C^{\top}_j
v\big)_{L^2(\partial\Omega)}\big], \quad  u,v\in H^{2m}(\Omega). 
\end{equation}

Next, following \cite{Gr68} and \cite{LM72}, we introduce the
spaces $D^s_A(\Omega) = \{u\in H^s(\Omega) \,|\,  Au\in H^0(\Omega)\}$,
$s\in \bbR$, provided with the graph norm  $\|u\|_{D^s_A(\Omega)}= (\|u\|^2_s
+ \|A u\|^2_0)^{1/2}$. Clearly,
$D^0_A(\Omega) = \dom(A)$ and $D^s_A(\Omega) \hookrightarrow H^s(\Omega)$.

\begin{definition}\label{defellip}
$(i)$ The operator ${\hatt A}_B$ defined by \eqref{2.54} is called
{\it elliptic} and is put in the class $\Ell(A)$ if $A$ is
properly elliptic on $\Omega$ and the system $\{B_j\}^{m-1}_{j=0}$ is
normal and  satisfies the covering condition (cf.\ \cite[Sects 2.1.1--2.1.4]{LM72})
at any point of the boundary $\partial\Omega$. \\
$(ii)$ The operator ${\hatt A}_B$ is called {\it coercive} in
$H^{s}(\Omega)$ with $s\ge2m$ if the {\it a priori} estimate (2.25) in \cite{Ag97} (cf.\ also \cite[Sect.\ 2.9.6]{LM72}) holds.
\end{definition}
We note that ${\hatt A}_B$ is a closed operator if $B$ satisfies the
covering condition (cf.\ \cite[Sect.\ 6.5]{Ag97}, \cite[Thm.\ 2.8.4]{LM72}).

If $\Omega$ is bounded, then any elliptic differential expression $\cA$ with
$C({\overline{\Omega}})$-coefficients is uniformly elliptic in  ${\overline{\Omega}}$. In this case, ${\hatt A}_B \in \Ell(A)$ if and only if ${\hatt A}_B$ is coercive in
$H^{2m}(\Omega)$ (see \cite{ADN59}, \cite[Sect.\ 2.9.6]{LM72}). If $\Omega$ is
unbounded, then the condition  ${\hatt A}_B \in \Ell(A)$ is still
necessary for coerciveness in $H^{2m}(\Omega)$, though, it is no
longer sufficient without additional assumptions on $A$.

\begin{hypothesis} \lb{h2.5} Assume that $\cA $ is a uniformly elliptic
operator, $a_{\alpha, \beta}(\cdot)\in C_b(\Omega)$ for $|\alpha| +
|\beta| \le 2m$ and $a_{\alpha, \beta}(\cdot)\in C_{ub}(\Omega)$ for
$|\alpha| + |\beta| = 2m$.
\end{hypothesis}

\begin{proposition}\label{prop4.2}
Assume Hypothesis \ref{h2.5},  ${\hatt A}_B\in \Ell(A)$,  and $0\in
\rho({\hatt A}_B)$. Then for any $s\in\bbR$, the mappings $B$
and $B^{\top}$ isomorphically map
$Z^s_A(\Omega) = \{u\in D^s_A(\Omega) \,|\, A_{\max}u=0\}$
and $Z^s_{A^{\top}}(\Omega) = \{u\in D^s_{A^{\top}}(\Omega)\,|\,  A^{\top}_{\max}u=0\}$
isomorphically onto $\Pi^m_{j=1} H^{s-m_j-(1/2)}(\partial\Omega)$ and
onto $\Pi^m_{j=1} H^{s-2m+\mu_j+(1/2)}(\partial\Omega)$, respectively.
\end{proposition}

\begin{definition} (\cite{Gr68, Vi63}) \label{def3.2A}
$(i)$ Under the assumptions of Proposition \ref{prop4.2}, let
$\varphi\in\Pi^{m-1}_{j=0} H^{s-m_j-(1/2)}(\partial\Omega)$, $s\in\bbR$.
Then one defines $P (z)\varphi$ to be the unique $u\in Z^s_{A-z I_{L^2(\Om)}}(\Omega)$
satisfying $B u=\varphi$. \\
$(ii)$ The {\it Calderon} operator $\Lambda (z)$ is defined  by
\begin{equation} \label{3.5A}
\Lambda (z) \colon  \Pi^{m-1}_{j=0} H^{s-m_j-(1/2)}(\partial\Omega) \to
\Pi^{m-1}_{j=0} H^{s-\mu_j-(1/2)}(\partial\Omega), \quad
\Lambda (z)\varphi = C P (z)\varphi.
\end{equation}
$(iii)$ Similarly, let $\psi\in\Pi^{m-1}_{j=0} H^{s-2m+\mu_j+1/2}(\partial\Omega)$. Then
$P (z)^\top \psi$ is defined
to be the unique solution in $Z^s_{A^{\top}-z I_{L^2(\partial\Om)}}(\Omega)$ of $B^{\top} u=\psi$ and
the {\it Calderon} operator $\Lambda (z)^\top$ is defined as
$\Lambda (z)^\top \psi=C^{\top} P^{\top}_z\psi$.
\end{definition}

Let $\Delta_{\partial\Omega}$ be the Laplace-Beltrami operator in
 $L^2(\partial \Omega)$, $-\Delta_{\partial\Om,1} =  -\Delta_{\partial\Omega}+I_{L^2(\partial\Om)}$. Then
 $-\Delta_{\partial\Om,1} = -\Delta_{\partial\Om,1}^* \geq I_{L^2(\partial\Om)}$. Moreover,
 $(-\Delta_{\partial\Om,1})^{-s/2}$ isomorphically maps $H^0 (\partial\Omega)$ onto
 $H^s (\partial\Omega)$, $s\in\bbR$.
Next, we introduce the diagonal $m\times m$ operator matrices
$-\Delta_{\partial\Om,1,m}$ and $-\Delta_{\partial\Om,1,\mu}$ with the  $(j,j)$-th entry
$(-\Delta_{\partial\Om,1})^{(m_{j}/2) +(1/4)}$ (resp., $(-\Delta_{\partial\Om,1})^{m-(\mu_{j}/2) -(1/4)})$.

\begin{proposition}\label{prop3.2}
Assume Hypothesis \ref{h2.5},  ${\hatt A}_B\in \Ell(A)$,  and $0\in \rho({\hatt A}_B)$.
Set
\begin{align}\label{3.10A}
\Gamma_{\Om,0} u & = (-\Delta_{\partial\Om,1,m})^{-1}Bu, \quad  \Gamma_{\Om,1} u =
(-\Delta_{\partial\Om,1,\mu}) (C u - \Lambda (0) Bu), \quad u\in \dom(A_{\max}), \\
\label{3.11A}
\Gamma^{\top}_{\Om,0} v& =(-\Delta_{\partial\Om,1,\mu})^{-1}B^{\top}v, \quad
\Gamma^{\top}_{\Om,1} v= \big(-\Delta_{\partial\Om,1,m}) (C^{\top}v-\Lambda (0)^\top B^{\top}v\big),
\quad v\in \dom(A^{\top}_{\max}).
\end{align}
Then the following holds: \\
$(i)$
$\Pi=\{\cH_{\partial\Om}\oplus\cH_{\partial\Om},\Gamma_{\Om},\Gamma^{\top}_{\Om}\}$,
with
$$
\cH_{\partial\Om} = \Pi_{j=0}^{m-1}H^0(\partial\Omega)=
\Pi_{j=0}^{m-1}L^2(\partial\Omega), \quad
\Gamma_{\Om}=(\Gamma_{\Om,0},\Gamma_{\Om,1}), \quad
\Gamma^{\top}_{\Om}=(\Gamma^{\top}_{\Om,0},\Gamma^{\top}_{\Om,1}),
$$
forms a boundary triple for the dual pair $\{A,A^{\top}\}$ of elliptic
operators in $L^2(\Om)$. In particular, the following Green's  formula holds
$$
(A_{\max}u,v)_{L^2(\Om)}-(u,A^{\top}_{\max}v)_{L^2(\Om)}
= (\Gamma_{\Om,1} u,\Gamma^{\top}_{\Om,0}
v)_{\cH_{\partial\Om}}-(\Gamma_{\Om,0} u,\Gamma^{\top}_{\Om,1} v)_{\cH_{\partial\Om}},
\quad u\in \dom(A_{\max}), \, v\in \dom(A^{\top}_{\max}).
$$
$(ii)$ The corresponding operator-valued Weyl--Titchmarsh function is given by
$$
M_{\Om,\Pi}(z)=(-\Delta_{\partial\Om,1,\mu})\bigl(\Lambda (z)-\Lambda (0)\bigr)
(-\Delta_{\partial\Om,1,m}), \quad z\in\rho({\hatt A}_B).
$$
\end{proposition}

In the context of operator-valued Weyl--Titchmarsh functions and
elliptic partial differential operators we also refer to the
recent preprint \cite{BGW08} (and the references cited therein).

\begin{definition}\label{def3.2}
For any operator $K \colon \dom(K) \to \Pi_{j=0}^{m-1} H^{-\mu_j - (1/2)}(\partial\Omega)$,
$\dom(K) \subseteq \Pi_{j=0}^{m-1} H^{-m_j - (1/2)}(\partial\Omega)$, we set
\begin{equation}\label{3.22B}
A_K  = A_{\max}\upharpoonright \dom(A_K ), \quad
\dom(A_K ) = \{u\in \dom(A_{\max}) \,|\,  Bu\in \dom(K), \  C u
= K B u\}.
\end{equation}
\end{definition}

\begin{definition}\label{def4.1}
Define
$\cS_p(\gH)=\{T\in \cB_\infty(\gH)\,|\, s_j(T) =O(j^{-1/p}) \, \text{as $j\to \infty$}\}$,
$p>0$, where $s_j(T)$, $j\in\bbN$, denote the singular values of $T$
(i.e., the eigenvalues of $(T^*T)^{1/2}$ ordered in decreasing magnitude, counting multiplicity).
\end{definition}

\begin{theorem}\label{th4.1A}
Assume the conditions of Proposition \ref{prop3.2} and suppose that
$0 \in \rho({\hatt A}_C)$ and $K\in \cC(\cH_{\partial\Om})$. Then: \\
$(i)$ For any realization $A_K \in \cC(L^2(\Om))$ of the form
\eqref{3.22B}, satisfying $\rho(A_K )\cap\rho({\hatt A}_B)\not =
\emptyset$, the following holds,
\begin{equation}\label{3.90}
\big[(A_{K} - z I_{L^2(\Om)})^{-l} - ({\hatt A}_B - z I_{L^2(\Om)})^{-l}\big] \in
\cS_{\frac{n-1}{2ml - 1/2}}\big(L^2(\Om)\big), \quad z\in\rho(A_{K})\cap\rho({\hatt A}_B),
\quad \ell\in\bbN.
\end{equation}
$(ii)$ If $B=\{B_j\}^{m-1}_{j=0}$ is a Dirichlet system, $K \in \cB(\cH_{\partial\Om})$,
and $\rho(A_K)\cap\rho({\hatt A}_B)\not =\emptyset$, one has
$$
\big[(A_K - z I_{L^2(\Om)})^{-1} - ({\hatt A}_B -z I_{L^2(\Om)})^{-1}\big] \in
\cS_{\frac{n-1}{2m}}\big(L^2(\Om)\big),\quad z\in\rho(A_K)\cap\rho({\hatt A}_B).
$$
\end{theorem}

Combining Weyl's theorem with Theorem \ref{th4.1A} one obtains the following result:

\begin{corollary}\label{cor4.4}
Assume the conditions of Theorem  \ref{th4.1A}. Then,
$\sigma_{\ess} (A_{K})=\sigma_{\ess}({\hatt A}_B)$.
\end{corollary}

In the case of elliptic realizations ${\hatt A}_G\in\Ell(A)$, we have the following stronger result:

\begin{theorem}\label{th3.smoothextens}
Suppose that the conditions of Proposition \ref{prop3.2}  are
satisfied and ${\hatt A}_{G}\in \Ell(A)$, that is, ${\hatt A}_{G}$ is
the elliptic realization of $A$ with $G=\{G_j\}_{j=0}^{m-1}$. Then for any $\ell\in \N$,
\begin{equation}\label{3.107B}
\big[({\hatt A}_{G} - z I_{L^2(\Om)})^{-\ell} - ({\hatt A}_B -z I_{L^2(\Om)})^{-\ell}\big] \in
\cS_{\frac{n-1}{2m\ell}}(L^2(\Om)), \quad z\in\rho({\hatt A}_{G})\cap\rho({\hatt A}_B).
\end{equation}
\end{theorem}

\section{The formally self-adjoint  case, nonnegative elliptic operators, and eigenvalues
in gaps}  \lb{s3}

Let $\cA$ be a formally self-adjoint  elliptic differential
expression of the form \eqref{2.10},  that is, $\cA =  \cA^{\top}$
or equivalently, $a_{p,q}= \overline{{a}_{q,p}} \in
C^{\infty}({\overline\Omega})$. In this case $A= A_{\min}=
A^{\top}_{\min} =A^{\top}$, that is, $A$ is symmetric, and
$A_{\max}= (A^{\top}_{\min})^*= A^*$. If a normal system
$\{B_j\}^{m-1}_{j=0}$ is chosen to be formally self-adjoint, that
is, $\hat{A}_B = (\hatt {A}_B)^*$,  then a system
$\{C_j\}^{m-1}_{j=0}$ can be chosen formally self-adjoint  too. In
this case $B_j^{\top}= B_j$ and $C_j^{\top}= C_j$,  and the
Green's  formula \eqref{2.19} holds with $B_j^{\top}= B_j$ and
$C_j^{\top}= C_j$. Moreover, in this case, $ \mu_j = \ord (C_j) = \ord (C^{\top}_j) = 2m-1-m_j$.
It follows that $\Delta_{\partial\Om,1,\mu} = \Delta_{\partial\Om,1,m}$. Hence, Proposition \ref{prop3.2} yields the following result:

\begin{proposition}\label{prop5.1}
Let $\cA$  be a formally symmetric elliptic differential expression and assume that
$\hatt {A}_B$ and $\hatt {A}_C$ are  self-adjoint. In addition,  assume the conditions of Propositon \ref{prop3.2} are satisfied. Then: \\
$(i)$ $\Pi=\{\cH_{\partial\Om},\Gamma_0,\Gamma_1\}$ with $ \cH_{\partial\Om} =
\Pi_{j=0}^{m-1}L^2(\partial\Omega)$,  and $\Gamma_{\Om,0},\, \Gamma_{\Om,1}$ defined
by \eqref{3.10A}, forms a boundary triple for the operator  $A^*$. In particular,
the following Green's  formula holds
$$
(A_{\max}u,v)_{L^2(\Om)} - (u,A_{\max}v)_{L^2(\Om)}
= (\Gamma_{\Om,1} u,\Gamma_{\Om,0}
v)_{\cH_{\partial\Om}}-(\Gamma_{\Om,0} u,\Gamma_{\Om,1} v)_{\cH_{\partial\Om}}, \quad u, v \in
D(A_{\max}).
$$
$(ii)$ The corresponding Weyl--Titchmarsh operator is given by
$M_{\Om,\Pi}(z)=(-\Delta_{\partial\Om,1,m})\bigl(\Lambda (z)-\Lambda (0)\bigr)
(-\Delta_{\partial\Om,1,m})$.
\end{proposition}

For any self-adjoint operator $T=T^*\in{\cC}(\gH)$ with associated family of spectral projections $E_T(\cdot)$, we set
$\kappa_{(\alpha,\beta)}(T) = \dim(E_T((\alpha,\beta))\gH)$, $-\infty \leq \alpha < \beta$ (these numbers may of course be infinite).

\begin{theorem}\label{th5.1}
Suppose that  $A> 0$ is a positive definite elliptic operator, and
$\Pi=\{\cH_{\partial\Om},\Gamma_0,\Gamma_1\}$ is the boundary
triple for $A^*$ in Proposition \ref{prop5.1} with $A_0 :=
A^*\upharpoonright \ker(\Gamma_0)  = {\hatt A}_{\gamma_D}$, the Dirichlet
realization of $A$. Assume also that the operator ${\hatt A}_C>0$
is positive definite, $0\in \rho({\hatt A}_C)$. Let  $K$ be a
densely defined $($not necessarily closed\,$)$ symmetric operator
in $\cH_{\partial\Om}$ and $A_K $  the corresponding extension
defined by \eqref{3.22B}. Then: \\
$(i)$ The Calderon operator  $\Lambda (0)$  is self-adjoint  and
negative definite, $\Lambda (0)<0$.  \\
$(ii)$ If  $K$ is  $\Lambda (0)$-bounded with bound strictly less than one, then
$A_K $ is symmetric $($but not  necessarily closed\,$)$.
\hspace*{5mm} If in addition, $\ran(K - \Lambda (0))$ is closed, then so is $A_K $,
that is, $A_K \in \cC(L^2(\Om))$.   \\
$(iii)$  If $K$ is  $\Lambda (0)$-compact and self-adjoint, then
$A_K $ is self-adjoint $A_K  = (A_K) ^*$, $\kappa_{(-\infty,0)}(A_K )<\infty $, and
\begin{equation}\label{4.36}
\kappa_{(-\infty,0)} (A_K )
= \kappa_{(-\infty,0)} \big(I_{L^2(\partial\Om)} + (-\Lambda (0))^{-1/2} K (-\Lambda (0))^{-1/2}\big).
\end{equation}
$(iv)$ If $K$ is  $\Lambda (0)$-compact and sectorial $($resp., $m$-sectorial\,$)$ with vertex $\zeta$ and semi-angle $\omega\in [0,\pi/2)$, then $A_K $ is
sectorial $($resp., $m$-sectorial\,$)$ with vertex $\zeta$ and semi-angle $\omega$ too.
\end{theorem}

\begin{proposition}\label{cor5.4}
Let  $\cA$ be formally self-adjoint  and assume the conditions of
Proposition \ref{prop5.1}. Assume also that $A_K  = (A_K) ^*$ is a
self-adjoint  extension of the form \eqref{3.22B} with $K\in
\cC(\cH)$. Then the absolutely continuous parts $A_{K, \ac}$ and
${\hatt A}_{B, \ac}$ of $A_K $ and ${\hatt A}_B$, respectively, are
unitarily equivalent. In particular, $\sigma_{\ac}(A_K ) =
\sigma_{\ac}({\hatt A}_B)$.
\end{proposition}

\begin{proposition}\label{cor5.5}
Suppose that $A=A_{\min}$ is  symmetric, and let ${\hatt A}_B =
{\hatt A}_{\gamma_D}$ be the Dirichlet realization of $A$. Assume
the conditions of Theorem \ref{th3.smoothextens} to be satisfied
and that ${\hatt A}_{G} = ({\hatt A}_{G})^* \in \Ell(A)$ is an
elliptic realization of $A$ with  $G= \{G_j\}_0^{m-1}$. Then  the
absolutely continuous parts ${\hatt A}_{G, \ac}$ and ${\hatt
A}_{\gamma_D,\ac}$ of ${\hatt A}_G$ and ${\hatt A}_{\gamma_D}$,
respectively, are unitarily equivalent. In particular,
$\sigma_{\ac}({\hatt A}_{G}) = \sigma_{\ac}({\hatt
A}_{\gamma_D})$.
\end{proposition}

Finally, we turn to eigenvalues in spectral gaps:

\begin{definition}
Let $A$ be a symmetric operator in $\cH$. Then
$(\alpha, \beta)$, $-\infty<\alpha < \beta < \infty$, is called a {\it gap of $A$} if
$\|(2A-(\alpha + \beta)I_{\cH})f\|_{\cH} \ge (\beta - \alpha)\|f\|_{\cH}$ for all $f\in\dom(A)$.
\end{definition}

By Corollary \ref{cor4.4}, $\sigma_\ess(A_K ) = \sigma_\ess({\hatt A}_B)$. Therefore,
in the gaps of ${\hatt A}_B$, the point spectrum of $A_K $ can possibly  accumulate  at most at the endpoints of the gaps. Next, we actually show that $\sigma_{\rm p}(A_K)$ cannot accumulate  at the left end point  of any  gap:

\begin{theorem}\label{th5.2}
Suppose that the conditions of Theorem \ref{th5.1}  are satisfied,
and that $K$ is a symmetric $\Lambda (0)$-compact operator in
$\cH_{\partial\Om}$. In addition, let $(\alpha, \beta)$ be  a
finite gap of $A_0= {\hatt A}_{\gamma_D}$
and introduce $T_0(z) = \Lambda (z)- \Lambda (0)$.  Then: \\
 $(i)$ $T(z) = \ol{T_0 (z)} \in \cB_\infty(\cH_{\partial\Om})$ for all
$z\in \rho({\hatt A}_{\gamma_D})$.   \\
$(ii)$ There exists $\varepsilon_0 \in (0, (\beta-\alpha)/2 )$
such that $E_{A_K}((\alpha, \alpha + \varepsilon_0)) =0$,  hence
$\kappa_{(\alpha, \beta- \varepsilon)} (A_K) = \dim(E_{A_K}
((\alpha, \beta-\varepsilon))) < \infty$ for any $\varepsilon\in
(0, \beta- \alpha)$.  Moreover, for any $\varepsilon\in (0,
\varepsilon_0)$ the following equality holds $($with $\Lambda
:=\Lambda (0)$$)$:
\begin{align*}
\kappa_{(\alpha, \beta- \varepsilon)} (A_K) & =
\kappa_{(-\infty,0)} (I_{\cH_{\partial\Om}} + \Lambda^{-1/2}(K +
T(\beta-
\varepsilon))\Lambda^{-1/2})   
- \kappa_{(-\infty,0)} (I_{\cH_{\partial\Om}} + \Lambda ^{-1/2} (K
+ T(\alpha + \varepsilon)) \Lambda ^{-1/2}).
\end{align*}
\end{theorem}

\begin{remark}
For Robin-type realizations $[\partial u/\partial n-\sigma
u]\upharpoonright\partial\Omega=0$, $\sigma \in
L^\infty(\partial\Om)$, of Schr\"odinger operators $-\Delta + q$
on exterior domains $\Omega\subset \bbR^3$, the estimate
\eqref{3.107B} (with $\ell=1$) goes back to the pioneering work by
Povzner \cite{Po67}. For Robin realizations $A_\sigma$ of a
second-order elliptic operator ${\mathcal A}= - \sum_{j,k=1}^n
\frac{\partial}{\partial x_j} a_{j,k}(x) \frac{\partial}{\partial
x_j} + q(x)$, with $q\geq 1$, and
$\sum_{j,k=1}^n \xi_j a_{j,k}(x) \xi_k >0$ for all $\{x, \xi\}\in
\Om \times (\bbR^n\setminus \{0\})$,
$$
A_\sigma = A_{\max} \upharpoonright\dom(A_\sigma),\; \dom(A_\sigma) =
\{u\in H^2(\Omega) \,|\,
({\partial u}/{\partial\nu} - \sigma u)\upharpoonright\partial\Omega =0 \}, \;
 {\partial}/{\partial\nu} = \sum_{j,k=1}^n a_{j,k}(x) \cos(n,
x_j) \f{\partial}{\partial x_k},
$$
$\sigma \in L^{\infty}(\partial\Omega)$, the estimate
\eqref{3.107B} was obtained by  Birman \cite{Bi62}. Moreover, in
\cite[Thm.\ 6.6]{Bi62} it is also proved that 
$\kappa_{(-\infty,0)} (A_\sigma) < \infty$. Thus, for $m=1$ and
$A_K  = A_\sigma$, equality \eqref{4.36} with $K$ being a
multiplication operator, $K \colon u\mapsto \sigma u$, yields a
stronger result as it describes the actual number of eigenvalues
in the gap $(-\infty, 0)$.

For positive elliptic realizations ${\hatt A}_G$ of a nonnegative
elliptic operator $A$ of order $2m$ in a bounded domain
$\Omega\subset \bbR^n$, the estimate \eqref{3.107B} is implied by
a sharp estimate due to Grubb \cite[eq.\ (3.22)]{Gr83}.
\end{remark}

Detailed proofs of these results will appear in \cite{GM09}.

\medskip

\noindent {\bf Acknowledgments.}
Mark Malamud is indebted to the Department of Mathematics of the University of Missouri, Columbia, MO, USA, for the hospitality extended to him during his month long stay in April/May of 2007 in connection with a Miller Family Scholarship, which enabled this collaboration.


\end{document}